\documentclass[A4paper]{article}
\usepackage{amsmath,amssymb,mathrsfs}
\date{}

\renewcommand{\mathcal}{\mathscr}
\def\cl#1{{\mathscr #1}}
\newcommand{\dlines}{\displaylines}
\newcommand{\field}[1]{\mathbb{#1}}

\newcommand{\finedim}{\par\hfill$\blacksquare$\par\noindent\ignorespaces}
\renewcommand{\th}{\theta}

\newcommand{\N}{\field{N}}

\newcommand{\Z}{\field{Z}}

\newcommand{\C}{\field{C}}

\newcommand{\e}{{\rm e}}

\def\tin#1{\par\noindent\hskip3em\llap{#1\enspace}\ignorespaces}
\def\cl#1{{\mathcal #1}}

\def\E{{\rm E}}

\def\P{{\rm P}}

\newtheorem{theorem}{Theorem}
\newtheorem{rema}[theorem]{Remark}

\newtheorem{assum}[theorem]{Assumption}

\newtheorem{cor}[theorem]{Corollary}

\newtheorem{example}[theorem]{Example}

\newtheorem{prop}[theorem]{Proposition}
\newtheorem{remark}[theorem]{Remark}
\newcommand{\proof}{{\it Proof}.\ \/}
\begin{document}
\title{\Huge\sc On the characterization of isotropic
Gaussian fields on homogeneous spaces of compact groups}
\author{P.Baldi, D.Marinucci\\
{\sl Dipartimento di Matematica, Universit\`a
di Roma {\it Tor Vergata}, Italy}\\
V.S.Varadarajan\\
{\sl Department of Mathematics, University of California at Los Angeles}
}
\maketitle

\begin{abstract}
\noindent Let $T$ be a random field invariant under the action of a compact group $G$
We give conditions ensuring that
independence of the random Fourier
coefficients is equivalent to Gaussianity. As a
consequence, in general it is not possible to simulate a
non-Gaussian invariant random field through its Fourier expansion using
independent coefficients.
\end{abstract}

\noindent{\it Key words and phrases} Isotropic Random Fields,
Fourier expansions, Characterization of Gaussian Random Fields.
\smallskip

\noindent{\it AMS 2000 subject classification:} Primary 60B15;
secondary 60E05,43A30.

\section{Introduction}\label{intro}%
Recently an increasing interest has been attracted by the topic of
rotationally  {\it real} invariant random fields on the sphere $\mathbb{S^2}$,
due to applications
to the statistical analysis of Cosmological and Astrophysical data
(see \cite{MR2065205}, \cite{mari2006a} and \cite{Kim}).

Some results concerning their structure and their spectral
decomposition have been obtained in \cite{BM06}, where a peculiar
feature has been pointed out, namely that if the
development into spherical harmonics
$$
T=\sum_{\ell=1}^\infty \sum_{-m}^m a_{\ell,m}Y_{\ell,m}
$$
of a rotationally invariant random field $T$ is such that the coefficients $a_{\ell,m}$, $\ell=1,2,\dots, 0\le
m\le \ell$ are independent, then the field is necessarily Gaussian
(the other coefficients are constrained by the condition
$a_{\ell,-m}=(-1)^m\overline a_{\ell,m}$). This fact (independence
of the coefficients+isotropy$\Rightarrow$Gaussianity) is not true
for isotropic random fields on other structures, as the torus or
$\Z$ (which are situations on which the action is Abelian).

This property implies in particular that non Gaussian rotationally
invariant random fields on the sphere {\it cannot} be simulated
using independent coefficients.

In this note we show that this
is a typical phenomenon for homogeneous spaces of compact
non-Abelian groups.
This should be intended as a contribution to a much
more complicated issue, i.e. the characterization of the isotropy of
a random field in terms of its random Fourier
expansion.

In \S2 and 3 we review some background material on harmonic analysis and
spectral representations for random fields. \S4 contains the main results, whereas
we moved to \S5 an auxiliary proposition.
\section{The Peter-Weyl decomposition}\label{Peter-Weil}%
Let $\mathcal X$ be a compact topological space and $G$ a
compact group acting on $\mathcal X$ transitively.  We
denote by $m_G$ the Haar measure of $G$. We know that there exists on
$\mathcal X$ a probability measure $m$ that is invariant by the
action of $G$, noted $x\to g^{-1}x$, $g\in G$. We assume that both
$m$ and $m_G$ are normalized and have total mass equal to $1$. We
shall write $L^2(\mathcal X)$ or simply $L^2$ instead of
$L^2(\mathcal X,m)$. Unless otherwise stated the spaces $L^2$ are
spaces of {\it complex valued} square integrable functions. We denote by
$L_g$ the action of $G$ on $L^2$, that is $L_gf(x)= f(g^{-1}x)$.

Let ${\widehat {\cl X}}$ be the set of equivalence classes of irreducible
unitary representations of $G$ which occur in the decomposition of
$L^2(\cl X, m)$. Since the action of $G$ commutes with the complex
conjugation on $L^2(\cl X, m)$, it is clear that for any irreducible
subspace $H$, $\overline H$, its conjugate subspace is also irreducible. If $H=\overline H$, we can find orthonormal bases $(\phi_k)$ for $H$ which are stable under conjugation; for instance we can choose the $\phi_k$ to be real. If $H\not=\overline H$, then there are two cases according as the action of $G$ on $\overline H$ is, or is not, equivalent to the action on $H$. If the two actions are inequivalent, then automatically $H\perp\overline H$. If the actions are equivalent, it is possible that $H$ and $\overline H$ are not orthogonal to each other. In this case $H\cap \overline H=0$ as both are irreducible and $S=H+\overline H$ is stable under $G$ and conjugation. In this case we can find $K\subset S$ stable under $G$ and irreducible such that $\overline K\perp K$ and $S=K\oplus \overline K$ is an orthogonal direct sum. The proof of this is postponed to the Appendix so as not to interrupt the main flow of the argument. We thus obtain the following orthogonal decomposition of $L^2(\cl X, m)$, {\it compatible with complex conjugation}:
\begin{equation}\label{e.pw-dec2}
L^2(\cl X, m)=\bigoplus_{i\in {\cl I}^o}H_i\oplus
\bigoplus_{i\in {\cl I}^+}(H_i\oplus \overline {H_i})
\end{equation}
where the direct sums are orthogonal and
$$
i\in {\cl I}^o\Leftrightarrow H_i=\overline H_i,\qquad
i\in {\cl I}^+\Leftrightarrow H_i\perp\overline H_i.
$$
We can therefore choose an orthonormal basis $(\phi_{ik})$ for $L^2(\cl X, m)$ such
that for $i\in {\cl I}^o$, $(\phi_{ik})_{1\le k\le d_i}$ is an orthonormal basis of
$H_i$ stable under conjugation, while, for $i\in {\cl I}^+$, $(\phi_{ik})_{1\le k\le d_i}$
is an orthonormal basis for $H_i$, where $d_i$ is the dimension of $H_i$; then, for $i\in {\cl I}^+$,
$(\overline {\phi_{ik}})_{1\le k\le d_i}$ is an orthonormal basis for $\overline H_i$. Such a orthonormal
basis $(\phi_{ik})_{ik}$ of $L^2(\cl X, m)$ is said to be {\it compatible with complex conjugation}.
\begin{example} \rm $\mathcal X=\mathbb{S}^1$, the one dimensional torus.
Here $\widehat G=\Z$ and $H_k$, $k\in\Z$ is generated by the
function $\gamma_k(\th)=\e^{ik\th}$. $\overline H_k=H_{-k}$ and
$\overline H_k\perp H_{k}$ for $k\not=0$. All of the $H_k$'s are
one-dimensional.
\end{example}
Recall that the irreducible representations of a compact topological group $G$ are all
one-dimensional if and only if $G$ is Abelian.
\begin{example} \rm  $G=SO(3)$, $\mathcal X=\mathbb{S}^2$, the sphere.
A popular choice of a basis of $L^2(\mathcal X,m)$ are the spherical
harmonics, $(Y_{\ell,m})_{-\ell\le m\le\ell}$, $\ell\in\N$ (see
\cite{MR1143783}). $H_\ell={\rm span}((Y_{\ell,m})_{-\ell\le
m\le\ell})$ are subspaces of $L^2(\mathcal X,m)$ on which $G$ acts
irreducibly. We have $\overline Y_{\ell,m}=(-1)^mY_{\ell,-m}$ and
$Y_{\ell,0}$ is real.

By choosing $\phi_{\ell,m}=Y_{\ell,m}$ for
$m\ge0$ and $\phi_{\ell,m}=(-1)^mY_{\ell,m}$ for $m<0$, we find a basis of
$H_\ell$ such that if $\phi$ is an element of the basis, then the
same is true for $\overline\phi$. Here $\dim(H_\ell)=2\ell+1$,
$\overline H_\ell=H_\ell$, so that in the decomposition
(\ref{e.pw-dec2}) there are no subspaces of the form $H_i$ for ${\cl I}^+$.
\end{example}
\section{The Karhunen-Lo\`eve expansion}
We consider on $\mathcal X$ a real {\it centered} square integrable
random field $(T(x))_{x\in \mathcal X}$. We assume that there exists
a probability space $(\Omega,\cl F,P)$ on which the r.v.'s $T(x)$ are
defined and that $(x,\omega)\to
T(x,\omega)$ is $\cl B(\cl X)\otimes\cl F$ measurable, $\cl B(\cl
X)$ denoting the Borel $\sigma$-field of $\cl X$. We assume
that
\begin{equation}\label{eq-l2bound}
\E\Bigl[\int_{\mathcal X}T(x)^2\, dm(x)\Bigr]=M<+\infty
\end{equation}
which in particular entails that $x\to T_x(\omega)$ belongs to $L^2(m)$ a.s.
Let us recall the main elementary facts concerning the Karhunen-Lo\`eve expansion
for such fields. We can associate to $T$ the bilinear form
on $L^2(m)$
\begin{equation}\label{eq-bilinear}
T(f,g)= \E\Bigl[\int_{\cl X}T(x) f(x)\, dm(x)
\int_{\cl X}T(y) g(y)\, dm(y)\Bigr]
\end{equation}
By (\ref{eq-l2bound}) and the Schwartz inequality one gets easily that
$$
|T(f,g)|\le M\Vert f\Vert_2\Vert g\Vert_2\ .
$$
Therefore, by the Riesz representation theorem there exists a function
$R\in L^2(\cl X\times\cl X,m\otimes m)$ such that
$$
T(f,g)=\int_{\cl X\times\cl X}f(x)g(y)R(x,y)\, dm(x)dm(y)\ .
$$
We can therefore define a continuous linear operator $R:L^2(m)\to L^2(m)$
$$
Rf(x)=\int_{\cl X} R(x,y)f(y)\,dm(y)\ .
$$
It can be even be proved that the linear operator $R$ is of trace class
and therefore compact (see \cite{MR2169627} for details). Since it is self-adjoint
there exists an orthonormal basis of $L^2(\cl X,m)$ that is formed by
eigenvectors of $R$.

Let us define, for $\phi\in L^2(\cl X,m)$,
$$
a(\phi)=\int_{\mathcal X}T(x)\phi(x)\, dm(x)\ ,
$$
Let $\lambda$ be an eigenvalue of $R$ and denote by $E_\lambda$ the
corresponding eigenspace. Then the following is well-known.
\begin{prop}\label{prop1}\sl Let $\phi\in E_\lambda$.
\tin{a)} If $\psi\in L^2(\cl X,m)$ is orthogonal to $\phi$, $a(\psi)$ is orthogonal
to $a(\phi)$ in $L^2(\Omega,\P)$. Moreover $\E[|a(\psi)|^2]=\lambda\Vert\psi\Vert_2^2$.
\tin{b)} If $\phi$ is orthogonal to $\overline\phi$,
then the r.v.'s $\Re a(\phi)$ and $\Im a(\phi)$ are orthogonal and have the
same variance.
\tin{c)} If the field $T$ is Gaussian, $a(\phi)$ is a Gaussian r.v. If
moreover
$\phi$ is orthogonal to $\overline\phi$,
then $a(\phi)$ is a complex centered Gaussian r.v.
(that is $\Re a_i$ and $\Im a_i$ are centered,
Gaussian, independent and have the same variance).
\end{prop}
\proof a) We have
$$
\dlines{
\E[a(\phi)\overline a(\psi)]=\E\Bigl[
\int_{\mathcal X}T(x)\phi(x)\, dm(x)
\int_{\mathcal X}T(y)\overline\psi(y)\, dm(y)\Bigr]=\cr
=\int_{\mathcal X\times \mathcal X}R(x,y)\phi(x)\overline\psi(y)
\, dm(x)\, dm(y)=\lambda\int_{\mathcal X}\phi(y)\overline\psi(y)
\, dm(y)=\lambda\langle \phi,\psi\rangle\ .\cr
}
$$
From this relation, by choosing first $\psi$ orthogonal to $\phi$ and then
$\psi=\phi$, the statement follows.
\tin{b)} From the computation in a), as $a(\overline\phi)=\overline{ a(\phi)}$, one gets
$\E[a(\phi)^2]=\lambda\langle \phi,\overline\phi\rangle$. Therefore, if
$\phi$ is orthogonal to $\overline\phi$, $\E[a(\phi)^2]=0$ which is equivalent
to $\Re a(\phi)$ and $\Im a(\phi)$ being orthogonal and having the
same variance.
\tin{c)} It is immediate that $a(\phi)$ is Gaussian. If $\phi$ is orthogonal to
$\overline\phi$, $a(\phi)$ is a complex centered Gaussian r.v., thanks to b).
\finedim
If $(\phi_k)_k$ is an orthonormal basis that is formed by
eigenvectors of $R$,
then under the assumption (\ref{eq-l2bound}) it is
well-known that the following expansion holds
\begin{equation}\label{e.kr-dev2}
T(x)=\sum_{k=1}^\infty a(\phi_k)\phi_k(x)
\end{equation}
which is called the Karhunen-Lo\`eve expansion.
This is intended in the sense of
$L^2(\cl X,m)$ a.s. in $\omega$. Stronger assumptions (continuity in
square mean of $x\to T(x)$, e.g.) ensure also that the convergence takes
place in $L^2(\Omega,\P)$ for every $x$ (see \cite{MR838963}, p.210 e.g.)


More relevant properties are true if we assume in addition that the random
field is invariant by the action $G$. Recall that the field $T$ is said to be (weakly) {\it invariant} by the action of $G$
if, for $f_i,\dots f_m\in L^2(\cl X)$ the joint laws of
$(T(f_1),\dots,T(f_m))$ and  $(T(L_g f_1),\dots,T((L_g f_m))$
are equal for every $g\in G$. Here we write
$$
T(f)=\int_{\cl X}T(x)f(x)\, dm(x),\qquad f\in L^2(\cl X)\ .
$$
If, in addition, the field is assumed to be continuous in sqare mean, this
imples that for every
$x_1,\dots, x_m\in\mathcal X$, 
$(T(x_1),\dots,T(x_m))$ and  $(T(g^{-1}x_1),\dots,T(g^{-1}x_m))$,
have the same joint laws for every $g\in G$. If the field is invariant then it is
immediate that the covariance function $R$ enjoys the invariance
property
\begin{equation}\label{inv:R}
R(x,y)=R(g^{-1}x,g^{-1}y)\qquad \strut\hbox to 0pt{a.e. for every $g\in G$}
\end{equation}
which also reads as
\begin{equation}\label{kh-inv}
 L_g(Rf)=R(L_gf)\ .
\end{equation}
Then, thanks to (\ref{kh-inv}), it is clear that
$G$ acts on $E_\lambda$. Therefore $E_\lambda$ is the direct sum of
some of the $H_i$'s introduced in the previous section.
Moreover it is a finite direct sum, unless $\lambda=0$, as the eigenvalues
of a compact operator that are different from $0$ cannot have but a finite
multiplicity. It turns out
therefore that the basis $(\phi_{ik})_{ik}$ of $L^2$ introduced in
the previous section is always formed by eigenvectors of $R$.

Moreover, if some of the $H_i$'s are of dimension $>1$, some of the
eigenvalues of $R$ have necessarily a multiplicity that is strictly
larger than $1$. As pointed out in \S \ref{Peter-Weil}, this
phenomenon is related to the non commutativity of $G$. For more details on
the Karhunen-Lo\`eve expansion and group representations see \cite{PecPyc2005h}.

%
%

Remark that if the random field is isotropic and satisfies
(\ref{eq-l2bound}), then (\ref{e.kr-dev2}) follows
by the Peter-Weyl theorem. Actually (\ref{eq-l2bound}) entails that,
for almost every $\omega$, $x\to T(x)$ belongs to $L^2(\cl X,m)$.
\begin{rema}\rm
An important issue when dealing with isotropic random fields is simulation.
In this regard, a natural starting point is
the Karhunen-Lo\`eve
expansion: one can actually sample random r.v.'s $\alpha(\phi_k)$,
(centered and standardized) and write
\begin{equation}\label{rem-sim}
T_n(x)=\sum_{k=1}^n\sqrt{\lambda_k}\,\alpha(\phi_k) \phi_k
\end{equation}
where the sequence $(\lambda_k)_k$ is summable. The point of course is
what conditions, in addition to those already pointed out,
should be imposed in order that (\ref{rem-sim}) defines
an isotropic field.
In order to have a real field, it will be necessary that
\begin{equation}\label{conjugate-cond}
\alpha(\overline\phi_k)=\overline{\alpha(\phi_k)}
\end{equation}
Our main result (see next section) is that if the $\alpha(\phi_k)$'s are independent r.v.'s
(abiding nonetheless to condition (\ref{conjugate-cond})), then the coefficients,
and therefore the field itself are Gaussian.
\end{rema}
If $H_i\subset L^2(\mathcal X,m)$ is a subspace on which $G$ acts
irreducibly, then one can consider the random field
$$
T_{H_i}(x)=\sum a(\phi_j)\phi_j(x)
$$
where the $\phi_j$ are an orthonormal basis of $H_i$. As remarked
before, all functions in $H_i$ are eigenvectors of $R$ associated to
the same eigenvalue $\lambda$.

Putting together this fact with (\ref{e.kr-dev2}) and
(\ref{e.pw-dec2}) we obtain the decomposition
$$
T=\sum_{i\in \mathcal I^\circ}T_{H^\circ_i}+\sum_{i\in \mathcal
I^+}(T_{H^+_i}+T_{H^-_i})\ .
$$
\begin{example}\rm
Let $T$ be a centered random field satisfying assumption (\ref{eq-l2bound})
over the torus $\mathbb{T}$,
whose Karhunen-Lo\`eve expansion is
$$
T(\th)=\sum_{-\infty}^{+\infty} a_k\,e^{ik\th},\qquad \th\in\mathbb{T}\ .
$$
Then, if $T$ is invariant by the action of $\mathbb{T}$ itself, the fields
$(T(\th))_\th$ and $(T(\th+\th'))_\th$ are equi-distributed, which implies
that the two sequences of r.v.'s
\begin{equation}\label{torus-inv1}
(a_k)_{-\infty<k<+\infty}\qquad\mbox{and}\qquad
(e^{ik\th'}a_k)_{-\infty<k<+\infty}
\end{equation}
have the same finite distribution for every $\th'\in\mathbb{T}$. Actually
one can restrict the attention to the coefficients
$(a_k)_{0\le k<+\infty}$, as necessarily $a_{-k}=\overline a_k$.

Conversely it is clear that if the two sequences in
(\ref{torus-inv1}) have the same distribution for every
$\th'\in\mathbb{T}$, then the field is invariant.

Condition (\ref{torus-inv1}) implies in particular that, for every
$k,-\infty<k<+\infty, k\not=0$ the distribution of $a_k$ must be
invariant by rotation (i.e. by the multiplication of a complex
number of modulus $1$).


If one assumes moreover that the r.v.'s $a_k,$ are independent, then
every choice of a distribution for $a_k, 0< k<+\infty $ that is
rotationally invariant gives rise to a random field that is
invariant with respect to the action of $\mathbb{T}$.
\end{example}
\section{$\!$Independent coefficients and non-Abelian groups}%
In this section we prove our main results showing that, if the group
$G$ is non commutative and under some mild additional assumptions,
independence of the coefficients of the Fourier development implies
their Gaussianity and, therefore, also that the random field must be
Gaussian. We stress that we {\it do not} assume independence of the
real and imaginary parts of the random coefficients.
\begin{prop}\label{plus} \sl Let $\cl X$ be an homogeneous space of the
compact group $G$.
Let $H^+_i\subset L^2(\mathcal X,m)$ be a
subspace on which $G$ acts irreducibly, having a dimension $\ge 2$
and such that if $f\in H^+_i$ then $\overline f\not\in H^+_i$. Let
$(\phi_k)_k$ be an orthonormal basis of $H^+_i$ and consider the
random field
$$
T_{H^+_i}(x)=\sum_k a_k\phi_k(x)\ .
$$
for a family of r.v.'s $(a_k)_k\subset L^2(\Omega,\P)$.
Then, if the r.v.'s $a_i$ are independent, the field $T_{H^+_i}$ is
$G$-invariant if and only if the r.v.'s $(a_k)_k$ are jointly
Gaussian and $\E(|a_k|^2)=c$ (and therefore also the field $T_{H^+_i}$ is Gaussian).
\end{prop}
\proof Since $G$ acts irreducibly on $H^+_i$, we have
$$
\phi_k(g^{-1}x)=\sum_{\ell=1}^{d_i}D_{k,\ell}(g)\phi_\ell(x)\ ,
$$
$d_i$ being the dimension of $H^+_i$ and $D(g)$ being the
representative matrix of the action of $g\in G$. Therefore
$$
T(g^{-1}x)=\sum_{\ell=1}^{d_i} \tilde a_\ell \phi_\ell(x)
$$
where
$$
\tilde a_\ell=\sum_{k=1}^{d_i} D_{k,\ell}(g)a_k\ .
$$
If the field is $G$-invariant, then the vectors $(\tilde
a_\ell)_\ell$ have the same joint distribution as $(a_k)_k$ and in
particular the $(\tilde a_\ell)_\ell$ are independent. One can then
apply
the Skitovich-Darmois theorem below (see 
\cite{MR0346969} e.g.) as soon as it is proved that $g\in G$ can be
chosen so that $D_{k,\ell}(g)\not=0$ for every $k,\ell$. This will
follow from the considerations below, where it is proved that the
set $Z_{k,\ell}$ of the zeros of $D_{k,\ell}$ has measure zero.

Indeed, let $G_1$ be the image of $G$ in the representation space so
that $G_1$ is also a connected compact group, and is moreover a Lie
group since it is a closed subgroup of the unitary group ${\rm
U}(d_i)$. If the representation is non trivial, then $G_1\not=\{1\}$
and in fact has positive dimension, and the $D_{k,\ell}$ are really
functions on $G_1$. For any fixed $k, \ell$ the irreducibility of
the action of $G_1$ implies that $D_{k,\ell}$ is not identically
zero. Indeed, if this were not the case, we must have $(g\phi_\ell,
\phi_k)=0$ for all $g\in G_1$, so that the span of the $g\phi_\ell$
is orthogonal to $\phi_k$; this span, being $G_1$-invariant and
nonzero, must be the whole space by the irreducibility, and so we
have a contradiction.

Since $D_{k\ell}$ is a non zero analytic function on $G_1$, it
follows from standard results that $Z_{k\ell}$ has measure zero.
Hence $\bigcup _{k\ell}Z_{k\ell}$ has measure zero also, and so its
complement in $G_1$ is non empty. \finedim
We use the following  version of the Skitovich-Darmois theorem,
which was actually proved by S.~G.~Ghurye and I.~Olkin
\cite{MR0137201} (see also \cite{MR0346969}).
\begin{theorem} \sl Let $X_1,\dots ,X_r$ be mutually
independent random vectors in $R^{n}.$ If the linear statistics
\begin{equation*}
L_{1}=\sum_{j=1}^{r}A_{j}X_{j}, \qquad L_{2}=\sum_{j=1}^{r}B_{j}X_{j}\ ,
\end{equation*}%
are independent for some real nonsingular $n\times n$ matrices $A_{j},B_{j},
$ $j=1,\dots ,r,$ then each of the vectors $X_{1},\dots ,X_{r}$ is normally
distributed.
\end{theorem}
We now investigate the case of the random field $T_H$, when $H$ is a
subspace such that $\overline H=H$. In this case we can consider a
basis of the form $\phi_{-k},\dots,\phi_k$, $k\le \ell$, with
$\phi_{-k}=\overline \phi_k$. The basis may contain a real function
$\phi_0$, if $\dim H$ is odd. Let us assume that the random
coefficients $a_k$, $k\ge 0$ are independent. Recall that $a_{
-k}=\overline{a_{k}}$.

The argument can be implemented along the
same lines as in Proposition \ref{plus}.
More precisely, if $m_1\ge 0$, $m_2\ge 0$, the two complex r.v.'s
\begin{equation}\label{tilde1}
\begin{array}{c}
\displaystyle \widetilde a_{m_{1}}=\sum_{m=-\ell}^\ell D_{m,m_{1}}(g)a_{m}\\
\displaystyle \widetilde a_{m_{2}}=\sum_{m=-\ell}^\ell D_{m
,m_{2}}(g)a_{m}
\end{array}
\end{equation}
%
have the same joint distribution as $a_{m_{1}}$ and $a_{ m_{2}}$.
Therefore, if $m_1\not =m_2$, they are independent. Moreover
$a_{-m}=\overline{a_{m}}$, so that the previous relation can be
written
\begin{align*}
\widetilde a_{m_{1}}&=D_{0 m_{1}}(g)a_{0}+ \sum_{m=1}^\ell
\Bigl(D_{m, m_{1}}(g)a_{m}+D_{-m, m_{1}}(g)
\overline{a_{m}}\Bigr)\\
\widetilde a_{m_{2}}&=D_{0, m_{2}}(g)a_{0}+ \sum_{m=1}^\ell\Bigl(
D_{m, m_{2}}(g)a_{m}+D_{-m, m_{2}}(g) \overline{a_{m}}\Bigr)
\end{align*}
In order to apply the Skitovich-Darmois theorem, we must ensure that
$g\in G$ can be chosen so that the real linear applications
\begin{equation}\label{condition}
z\to D_{m, m_{i}}(g)z+D_{-m, m_i}(g)\overline z, \qquad
m=1,\dots,\ell, i=1,2
\end{equation}
are all non singular. It is immediate that this condition is
equivalent to imposing that $|D_{m, m_{i}}(g)|\not=|D_{-m,
m_{i}}(g)|$.

We show below that (\ref{condition}) is satisfied for some
well-known examples of groups and homogeneous spaces. We do not know
whether (\ref{condition}) is always satisfied for every compact
group. We are therefore stating our result conditional upon
(\ref{condition}) being fulfilled.
\begin{assum} \label{assum0}There exist $g\in G$, $0\le m_1<m_2\le\ell$ such that
$$
|D_{m, m_{i}}(g)|\not=|D_{-m, m_{i}}(g)|
$$
for every $0\le m\le \ell$.
\end{assum}
We have therefore proved the following.
\begin{prop}\label{zero} \sl Let $\cl X$ be an homogeneous space of the
compact
group $G$. Let $H_i\subset L^2(\mathcal X,m)$ be a
subspace
on which $G$ acts irreducibly, having a dimension $d> 2$ and such that
$\overline H_i=H_i$. Let $(\phi_k)_k$ be an
orthonormal basis of $H_i$ such that $\phi_{-k}=\overline \phi_k$
and consider the random field
$$
T_{H_i}(x)=\sum_k a_k\phi_k(x)
$$
where the r.v.'s $a_k,k\ge 0$ are centered, square integrable, independent and
$a_{-k}=\overline a_k$. Then $T_{H_i}$ is $G$-invariant if and
only if the r.v.'s $(a_k)_{k\ge 0}$ are jointly Gaussian and
$\E(|a_k|^2)=c$ (and therefore also the field $T_{H_i}$ is
Gaussian).
\end{prop}
%
Putting together Propositions \ref{plus} and \ref{zero} we obtain
our main result.
\begin{theorem}\label{main} Let $\cl X$ be an homogeneous space of the
compact group $G$. Consider the decomposition (\ref{e.pw-dec2}) and let
$\big((\phi_{ik})_{i\in \mathcal I^\circ}, (\phi_{ik},
\overline\phi_{ik})_{i\in \mathcal I^+} \big)$ be a basis of
$L^2(G)$ adapted to that decomposition. Let
$$
T=\sum_{i\in \mathcal I^\circ}\sum_k a_{ik}\phi_{ik}+\sum_{i\in
\mathcal I^+}\sum_k \big(a_{ik}\phi_{ik}+\overline a_{ik}\overline
\phi_{ik}\big)
$$
be a random field on $\cl X$, where the series above are intended to
be converging in square mean. Assume that $T$ is isotropic with
respect to the action of $G$ and that the coefficients $(a_{ik})_{i\in
\mathcal I^\circ, k\ge 0}, (a_{ik})_{i\in I^+}$ are independent.
If moreover
\tin{a)} the only one-dimensional irreducible representation appearing
in (\ref{e.pw-dec2}) are the constants;
\tin{b)} there are no $2$-dimensional subspaces $H\subset L^2(\cl X)$,
invariant under the action of $G$ and such that $\overline H=H$;
\tin {c)} The random coefficient corresponding to the trivial
representation vanishes.
\tin{d)} For every $H\subset L^2(\cl X)$, irreducible under the action of $G$ and such that
$\overline H=H$, Assumption \ref{assum0} holds.

Then the coefficients $(a_{ik})_{i\in \mathcal I^\circ, k\ge 0},
(a_{ik})_{i\in \mathcal I^+}$ are Gaussian and the field itself is
Gaussian.
\end{theorem}
Let us stress with the following statements the meaning of assumption a)--d).
The following result gives a condition ensuring that assumption b) of
Theorem \ref{main} is satisfied.
\begin{prop}\label{prop411} Let $U$ be an irreducible unitary $2$-dimensional
representation of $G$ and let $H_1$ and $H_2$ be the two
corresponding subspaces of $L^2(G)$ in the Peter-Weyl decomposition.
Then if $U$ has values in $SU(2)$, then
$\overline{H}_1=H_2\not=H_1$.
\end{prop}
\proof If we note
$$
U(g)=\begin{pmatrix}a(g)&b(g)\cr
c(g)&d(g)\cr
\end{pmatrix}
$$
then one can assume that $H_1$ is generated by the functions $a$ and
$c$, whereas $H_2$ is generated by  $b$ and $d$. It suffices now to
show that $\overline a$ is orthogonal both to $a$ and $c$. But, the
matrix $U(g)$ belonging to $SU(2)$, we have $\overline{a(g)}=d(g)\in
H_2$. \finedim
Recall that the commutator $G_0$, of a topological group $G$ is the closed
group that is generated by the elements of the form $xyx^{-1}y^{-1}$
\begin{cor}\label{semisimple} Let $G$ be a compact group
such that its commutator $G_0$ coincides with $G$ himself. Then assumptions
a) and b) of Theorem \ref{main} are satisfied. In particular these assumptions
are satisfied if $G$ is a semisimple Lie group.
\end{cor}
\proof Recall that if $G_0=G$, $G$ cannot have a
quotient that is an abelian group. If there was a unitary
representation $U$ with a determinant not identically equal to $1$,
then $g\to\det(U(G))$ would be an homomorphism onto the torus
$\mathbb{T}$ and therefore $G$ would possess $\mathbb{T}$ as a
quotient. The same argument proves that $G$ cannot have a one
dimensional unitary representation other than the trivial one.
One can therefore apply Proposition \ref{prop411} and b) is satisfied.
\finedim
\begin{remark} \rm It is easy to prove that Assumption \ref{assum0} is satisfied
when $\cl X=\mathbb S^2$ and $G=SO(3)$. As mentioned in \cite{BM06},
this can be established using explicit expressions of the representation
coefficients as provided e.g. in \cite{MR1022665}.

In the same line of arguments it is also easy to check the same in the cases $\cl X=SO(3)$, $G=SO(3)$ and
$\cl X=SU(2)$, $G=SU(2)$.
\end{remark}
\begin{remark} \rm As far as condition c) of Theorem \ref{main}, let us remark that
the coefficient of the trivial representation
corresponds to the empirical mean of the field. As any random field can be decomposed
into the sum of its empirical mean plus a field whose coefficient of the trivial representation
vanishes, our result can be interpreted in terms of Gaussianity of this second component.
\end{remark}
\section{Appendix}\label{sec-appendix}%
\begin{prop}\label{raja5}
Let $V$ be a finite dimensional Hilbert space on which $G$ acts
unitarily, and let $V$ be equipped with a conjugation $\sigma (v\to
\bar v)$ commuting with the action of $G$. Let $H$ be an irreducible
$G$-invariant subspace and let $V=H+\overline H$. \tin{a)} If the
actions of $G$ on $H$ and $\overline H$ are inequivalent, then
$\overline H\perp H$ and $V=H\oplus \overline H$. \tin{b)} If the
actions of $G$ on $H$ and $\overline H$ are equivalent, then either
$H=\overline H$ or we can find an irreducible $G$-invariant subspace
$K$ of $V$ such that $\overline K\perp K$ and $V=K\oplus \overline
K$.
\end{prop}
\proof Let $P$ be the orthogonal projection $V\to \overline H$ and
$A$ its restriction to $H$. Then, for every $h\in H$,
$h'\in\overline H$ and $g\in G$, we have
$$
\langle g(Ah),h'\rangle=\langle Ah,gh'\rangle=\langle
h,gh'\rangle=\langle gh,h'\rangle=\langle A(gh),h'\rangle
$$
From this we get that $G$ acts on $A(H)$. The action of $G$ on
$\overline H$ being irreducible, we have either $A(H)=\{0\}$ or
$A(H)=\overline H$. In the first case $H$ is already orthogonal to
$\overline H$. Otherwise $A$ intertwines the actions on $H$ and on
$\overline H$, so that these are equivalent and $V=H\oplus H^\perp$.

$V$ being the sum of two copies of the representation on $H$, there
is a {\it unitary} isomorphism $V\simeq H\otimes \C^2$ where $\C^2$
is given the standard scalar product. So we assume that $V=H\otimes
\C^2$. $G$ acts only on the first component, so that $G$ acts
irreducibly on every subspace of the form $H\otimes Z$, $Z$ being a
one dimensional subspace of $\C^2$.

Let us identify the action of $\sigma$ on $H\otimes \C^2$. Let
$\sigma_0$ be the conjugation on $V$ defined by $\sigma_0(u\otimes
v)=u\otimes \bar v$ where $v\to \bar v$ is the standard conjugation
$(z_1,z_2)\to (\overline {z_1}, \overline {z_2})$. Then
$\sigma\sigma_0$ is a {\it linear operator\/} commuting with $G$ and
so is of the form $1\otimes L$ where $L(\C^2\to \C^2)$ is a linear
operator. Hence
$$
\sigma (u\otimes v)=\sigma\sigma_0 (u\otimes \overline v)=u\otimes
L\bar v.
$$
If $Z$ is any one dimensional subspace of $\C^2$, $H\otimes Z$ is
$G$-invariant and irreducible, and we want to show that for some
$Z$, $H\otimes Z\perp H\otimes Z^\sigma$, i.e., $Z\perp Z^\sigma$.
Here $ Z^\sigma=\sigma(Z)$.

For any such $Z$, let $v$ be a nonzero vector in it; then the
condition $Z\perp Z^\sigma$ becomes $(v, L\bar v)=0$ where $(,)$ is
the scalar product in $\C^2$. Since $(,)$ is Hermitian,
$B(v,w):=(v,L\bar w)$ is {\it bilinear} and we want $v$ to satisfy
$B(v,v)=0$. This is actually standard: indeed, replacing $B$ by
$B+B^T$ (which just doubles the quadratic form) we may assume that
$B$ is symmetric.

If $B$ is degenerate, there is a nonzero $v$ such that $B(v,w)=0$
for all $w$, hence $B(v,v)=0$. If $B$ is nondegenerate, there is a
basis $v_1, v_2$ for $\C^2$ such that $B(v_i,v_j)=\delta_{ij}$.
Then, if $w=v_1+iv_2$ where $i=\sqrt {-1}$, $B(w,w)=0$.\finedim

\bibliography{bibbase}
\bibliographystyle{alpha}
\end{document}